\documentclass[10pt]{article}
\usepackage{amscd,amsfonts,amsmath,amssymb,amstext,amsthm}
\title{Cycles in $G$-orbits in $G^\mathbb C$-flag manifolds}

\date{ }
\author{A. Huckleberry\footnote{
Research partially supported by Schwerpunkt "Global methods in complex
geometry" and SFB-237 of the Deutsche Forschungsgemeinschaft.}
\,\, \& B. Ntatin\footnote{
Supported by a stipend of the Deutsche Akademische Austauschdienst.}
}

\newcommand{\s} {{\smallskip \noindent}}
\newcommand{\m} {{\medskip \noindent}}

\newtheorem{thm} {Theorem} 
\newtheorem{lem} [thm]{Lemma}
\newtheorem{prop}[thm]{Proposition}
\newtheorem{cor} [thm]{Corollary}
\begin{document}
\maketitle
\abstract{\begin{quote} \footnotesize
There is a natural duality between orbits $\gamma $ of a real form $G$ 
of a complex semisimple group $G^\mathbb C$ on a homogeneous 
rational manifold $Z=G^\mathbb C/P$ and those $\kappa $ of the
complexification $K^\mathbb C$ of any of its maximal compact
subgroups $K$: $(\gamma ,\kappa )$ is a dual pair if $\gamma \cap \kappa $
is a $K$-orbit. The cycle space $C(\gamma )$
is defined to be the connected component containing the identity
of the interior of $\{ g:g(\kappa )\cap \gamma  \  \text{is non-empty and
compact} \}$.  Using methods which were recently developed 
for the case of open $G$-orbits, geometric properties of 
cycles are proved, and it is shown that $C(\gamma )$ is contained in a 
domain defined by incidence geometry. In the non-Hermitian case 
this is a key ingredient for proving that $C(\gamma )$ is a certain explicitly 
computable universal domain

\end {quote}}
\setcounter{section}{-1}
\section {Introduction and Notation}

Let $G$ be a non-compact semi-simple Lie group without compact factors
which is embedded in its complexification $G^\mathbb C$ and
let $Z=G^\mathbb C/Q$ be a $G^\mathbb C$-flag manifold, i.e., a compact,
homogeneous, algebraic $G^\mathbb C$-manifold. Denote by $K$ a 
maximal compact subgroup of $G$; in particular $G/K$ is a 
negatively curved Riemannian symmetric space.  In the sequel
we will assume that $G$ is simple.  The necessary adjustments 
for the semi-simple case are straight-forward.

\m
Let $Orb_Z(G)$ (resp. $Orb_Z(K^\mathbb C)$) denote the set of $G$-orbits
(resp. $K^\mathbb C$-orbits) in $Z$. It is known that these
sets are finite (\cite {W1}).
If $\kappa \in Orb_Z(K^\mathbb C)$ and $\gamma \in Orb _Z(G)$, then
$(\kappa ,\gamma )$ is said to be a dual pair if $\kappa \cap \gamma $
is non-empty and compact.

\s
If $\gamma $ is an open $G$-orbit, then $\kappa $ being dual to 
$\gamma $ is equivalent to $\kappa \subset \gamma $. In (\cite {W1})
it is shown that every open $G$-orbit contains a unique compact
$K^\mathbb C$-orbit, i.e., duality at the level of open $G$-orbits.
This is extended in (\cite {M}, see also \cite {BL} and \cite {MUV}) 
to the case of all orbits: 
{\it For every $\gamma \in Orb_Z(G)$ there exists a unique 
$\kappa \in Orb_Z(K^\mathbb C)$ such that $(\gamma ,\kappa )$ 
is a dual pair and vice versa.} Furthermore, if $(\gamma ,\kappa )$ is
a dual pair then the intersection $\kappa \ \cap \ \gamma $ is transversal
at each of its points and consists of exactly one $K$-orbit.

\m
To motivate the notion of a cycle in this case, let us begin with the
case of an open $G$-orbit $D$.  The dual orbit $\kappa $
defines a point in the cycle space ${\cal C}^q(D)$,
where $q:=dim_\mathbb C\kappa $.  The connected component
$\Omega _W(D)$ of
$\{g\in G^\mathbb C:g(\kappa )\subset D\}$
can be regarded as a family of cycles by the procedure of associating
$g$ to the cycle $g(\kappa )$. Clearly $\Omega _W(D)$ is invariant
by the action of $K^\mathbb C$ on $G^\mathbb C$ on the right and therefore
we often regard it as being in the affine homogenous space
$\Omega :=G^\mathbb C/K^\mathbb C$.

\m
The cycle space $C(\gamma )$ associated to a lower-dimensional orbit
is defined analogously, at least when one has 
duality in mind:
$C(\gamma )$ is the connected component containing the identity of
the interior of the set 
\begin {equation*}
\{ g\in G^\mathbb C: g(\kappa )\cap \gamma \text { is non-empty and compact}\}. 
\end {equation*}

Since the intersection $\kappa \cap \gamma $ is transversal, it is
clear that $C(\gamma )$ is non-empty

\m
Here, in a result that we state at the outset, it is shown
that this is in fact a reasonable set.  The following is proved
in $\S\, \ref {transversality}$.

\begin {prop} \label {real cycles}
If $g\in C(\gamma )$, then
$g(\text{c}\ell (\kappa ))\cap \text{c}\ell (\gamma )=:M_g$ is
a compact subset of $g(\kappa )\cap \gamma $. This intersection
is transversal at each of its points and the manifold $M_g$ is
homotopic in $\gamma $ to $M=M_e$.  If $g\in \text{bd}(C(\gamma ))$, 
then either $g(\kappa )\cap \gamma $ is non-compact or empty. In 
particular, if $\{ g_n\} $ is a divergent sequence in $C(\gamma )$, then
there exists $z_n\in g_(\kappa )\cap \gamma $ with $\{ z_n\} $ divergent
in $\gamma $.
\end {prop}

In the sequel we take this to be the meaning of the {\it cycle space}
$C(\gamma )$ and regard the manifolds $M_g$ as cycles.  Just as in the case
of open orbits, we often regard $C(\gamma )$ as a $G$-invariant
domain in $\Omega $.  

\m
Our main result is the complete characterization of the 
cycle spaces $C(\gamma )$ in the case where $G$ is non-Hermitian, where
our contribution is in the case of the non-open orbits.
For a given group $G$ these cycle spaces are all naturally biholomorphically
isormorphic to a fixed domain which, when realized in $\Omega $,
is denoted by $\Omega _{AG}$ (see Thm. \ref {equality}).

\m
This domain was discussed in (\cite {C}) in the context of differential
geometry. If $G\times _K{\mathfrak p}$ is the tangent
bundle of the Riemannian symmetric space $M=G/K$, then $\Omega _{AG}$
is identified by the polar coordinates map ${\cal P}:TM\to \Omega $, 
$[(g,\xi)]\mapsto gexp(i\xi )$, with the maximal neighborhood
of the $0$-section on which ${\cal P}$ is a local diffeomorphism.
It was considered in (\cite {AG}) from the point of view of
neigborhoods of $M$ in $\Omega $ on which $G$ acts properly.  It is
known to be a Stein domain and its complex geometry is closely related
to the Riemannian geometry of the symmetric space $M$ (\cite {BHH}).

\s
In terms of roots, if ${\mathfrak g}={\mathfrak k}+{\mathfrak a}
+{\mathfrak n}$ is an 
Iwasawa-decomposition, $\Phi $ is a system of roots on ${\mathfrak a}$
and $\omega _{AG}$ is the connected component containing $0\in {\mathfrak a}$
of the set which is obtained by removing from 
${\mathfrak a}$ the union affine hyperplanes

\begin {equation*}
\underset{\alpha \in \Phi}{\cup }\{\xi \in {\mathfrak a}: 
\alpha (\xi )=\frac{\pi }{2}\}.
\end {equation*}  

Then $\Omega _{AG}=Gexp(i\omega _{AG}).x_0$,
where $x_0\in G^\mathbb C/K^\mathbb C$ is the base point.

\m
Using a certain triality and related incidence varieties (\cite {HW1}) 
along with $G$-invariant theoretic properties and the Kobayashi
hyperbolicity of $\Omega _{AG}$ (\cite {H},\cite {FH}), it has
been recently shown that, with a few well-known Hermitian exceptions
where $\Omega _W(D)$ is the associated bounded symmetric domain, the cycle
domain $\Omega _W(D)$ of an open orbit is naturally identifiable
with $\Omega _{AG}$. 

\s
Our work here makes use of the methods
and results of these papers along with a result of (\cite {GM}) which,
together with knowledge of the intersection of the cycle domains for
the open orbits in $G^\mathbb C/B$, implies the inclusion
$\Omega _{AG}\subset C(\gamma )$ for all $\gamma \in Orb_Z(G)$.

\s
For a survey of these and other basic properties of 
$\Omega _{AG}$ see (\cite {HW2}).   

\section {Basic triality} \label {triality}

Here $B$ denotes a Borel subgroup of $G^\mathbb C$ which
contains the factor $AN$ of an Iwasawa decomposition $G=KAN$.
As usual the closure $S$ of a $B$-orbit $O$ in $Z$ is referred to as
a Schubert variety.  Let $Y:=S\setminus O$.

\m
Since the set of Schubert varieties ${\cal S}$ generates the homology of $Z$,
for every $\kappa \in Orb_Z(K^\mathbb C)$ the set
\begin {equation*}
{\cal S}_\kappa :=\{ S\in {\cal S}:S\cap \text{c}\ell (\kappa )\not =
\emptyset \}
\end {equation*}
is non-empty. The following is a slightly refined version of Thm. 3.1
of (\cite {HW1}).

\begin {thm}
Let $(\kappa ,\gamma )$ be a dual pair and $S\in {\cal S}_\kappa $.
Then

\begin {enumerate}

\item $S\cap \text{c}\ell (\kappa )\subset \kappa \cap \gamma $
\item The intersection $S\cap \gamma $ is open in $S$ and consists
of finitely many $AN$-orbits each of which is open in $S$ and closed 
in $\gamma $.
\item If $z_0\in S\cap \gamma $, then $\Sigma :=(AN).z_0$ intersects
$\kappa  $ in exactly one point and that interesection is transversal
in $Z$.

\end {enumerate}
\end {thm}

The only refinement is that, instead of $\Sigma \cap \kappa $
being finite, we now show that it consists of just one point.
This is implicit in \cite {HW1} (see $\S 5$ in that paper). 
Let us repeat the relevant details.

\m
Let $\alpha :K_{z_0}\times (AN)_{z_0}\to G_{z_0}$ be defined
by multiplication, $(k,an)\mapsto kan$.  As was shown in
(\cite {HW1}) $\alpha $ is a diffeomorphism onto a number of
components of $G_{z_0}$.  However, since the orbit $M:=K.z_{0}$ 
is a strong deformation retract of $\gamma $, it follows that
$G_{z_0}/K_{z_0}$ is connected and consequently $\alpha $ is
surjective.  Thus $\vert \Sigma \cap \kappa \vert =1$ is
proved just as in Cor. 5.2 of (\cite {HW1}). $\square $

\m
Given $\gamma $,$\kappa $, an Iwasawa-decomposition $G=KAN$,
a Borel subgroup $B$ of $G^\mathbb C$ which contains $AN$,
a $B$ Schubert variety $S\in {\cal S}_\kappa $ and an
intersection point $z_0\in \kappa \cap S$ as above, we
refer to $\Sigma =AN.z_0$ as the associated {\it Schubert slice}.
Note that $G$-conjugation yields a $G_p$-invariant family 
of Schubert slices $\Sigma $ at every point $p\in \gamma $. 

\s
It should also be noted that in the case that $\gamma =
\gamma _{\text{c}\ell }$ is closed, $\Sigma $ is just 
the associated $B$-fixed point.

\section {Cycle transversality} \label {transversality}

Our main goal here is to prove Prop.$\, 1$. The first step
is to show that if $g(\kappa )\cap \gamma $ is compact, the
intersection of the variety $g(\text{c}\ell(\kappa ))$ with
every Schubert slice $\Sigma $ has optimal transversality properties.

\begin {lem} \label {intersection}
For all $g\in G^\mathbb C$ with $g(\kappa )\cap \gamma $ 
compact the number of
points in the intersection $g(\kappa )\ \cap  \ \Sigma $
is bounded by the intersection number $S.\text{c}\ell(\kappa )$.
\end {lem}
\begin {proof} Since $\Sigma $ can be regarded as a domain in
$O\cong \mathbb C^n$ and $g(\kappa )\cap \gamma $ is compact,
it follows from the maximum principle that $g(\kappa )\cap \Sigma $
is finite and of course it is then bounded by the the intersection
number $S.\text{c}\ell (\kappa )$.
\end {proof}

It should be underlined that, 
since $\vert g(\kappa )\cap \Sigma \vert $ is finite, it is semi-continuous
in the sense that it can only increase as $g$ moves away from $g_0$.

\m
Let
\begin {equation*}
I:=\{ g\in G^\mathbb C:\vert g(\kappa )\cap \Sigma \vert =1
\text{ for all } \Sigma \}
\end {equation*}

Here {\it for all \ $\Sigma $} means for all choices of the maximal
compact group $K$ and all Iwasawa factors $AN$, 
i.e., all Schubert slices which arise by
$G$-conjugation of those $\Sigma $ which are connected components
of $S \cap \gamma $ for a fixed $S\in {\cal S}_\kappa $. 

\m
We now consider the open set $I\cap C(\gamma )$. 

\begin {lem} \label {main lemma}
If $g\in \text{bd}(C(\gamma )\cap I)$, then either 
$g(\kappa )\cap \gamma $ is empty or non-compact.  In particular,
$C(\gamma )\subset I$.
\end {lem}
\begin {proof}
We may assume that $g\not \in I$, because if $g\in I$ and
$g(\kappa )\cap \gamma $ is compact, then both conditions hold
in an open neighborhood of $g$ and consequently $g\in C(\gamma )$.

\s
We assume that $g(\kappa )\cap \gamma $ is compact
and reach a contradiction.  For this let $\tilde M$ be the
limiting set of the sequence of manifolds $M_n:=g_n(\kappa )\cap \gamma $.
Since $C(\gamma )$ is by definition connected and $M_e$ is connected,
it follows that $\tilde M$ is a connected closed set.

\m
Since $g(\kappa )\cap \gamma $
is compact, it follows that $\tilde M\cap \gamma $ is compact and is 
therefore closed in $\tilde M$.  On the other hand 
$\tilde M\cap \text{bd}(\gamma )$ is
closed in $\tilde M$ and consequently $\tilde M\cap \gamma $
is also open in $\tilde M$.  By the semicontinuity of the
intersection number and the fact that $g\not \in I$, 
$\tilde M\cap \Sigma _0=\emptyset $ for some Schubert slice
$\Sigma _0$. In particular, since $\vert M_n\cap \Sigma _0\vert =1$ for all $n$,
$\tilde M\not \subset \gamma $.
Since $\tilde M$ is connected, it follows that 
$\tilde M \cap \gamma =\emptyset $.

\m
But $g(\kappa )\cap \gamma $ is non-empty.  So there
exists $p\in \gamma $ and an open neighborhood $U(p)$ so that
$g_n(\kappa )\cap U=\emptyset$, but $p\in g(\kappa )$.  On the
other hand, if $\Sigma $ is a Schubert slice at $p$, 
and $p$ were isolated in $g(\kappa )\cap \Sigma $, it would 
follow that $g_n(\kappa )\cap U\cap \Sigma $ is non-empty for 
$n$ sufficiently large.
\end {proof}

\m
{\it Proof of Proposition \ref {real cycles}}. Observe that 
if $g\in C(\gamma )$ and
$g(\text{c}\ell (\kappa ))\cap \text{c}\ell (\gamma )$ contained
an additional point $p\in \text{bd}(\gamma )$ or in $\text {bd}(\kappa )$,
then for $h\in C(\gamma )$ chosen appropriately, in particular small,
we would find some $\Sigma _0$ which contains $h(p)$ as well as
another point in $M_{hg}$. This is contrary to the fact that
$C(\gamma )\subset I$ and that this intersection contains exactly one point.  

\s
The transversality of the intersection and the properties of $M_g$
are also immediate consequences of $C(\gamma )$ being connected and
contained in $I$.

\s
If $g\in \text{bd}(C(\gamma ))$, then, again 
since $C(\gamma )\subset I$, by Lemma \ref {main lemma}
$g(\kappa )\cap \gamma $ is either empty or non-compact.
\hfill $\square $

\section {Description of the cycle spaces}

\subsection {Lifting cycle spaces}

The above transversality results are only useful in the
case when $\gamma $ is not closed. Here we prove a 
{\it Lifting Lemma} which will be used to understand
the cycle space $C(\gamma _{\text{c}\ell })$ of the
closed orbit.  For this let $Z=G^\mathbb C/P$ as usual, 
let $\tilde Z=G^\mathbb C/\tilde P$ be
defined by a parabolic group $\tilde P$ which is contained in
$P$ and $\pi :\tilde Z\to Z$ is the natural projection.

\begin {prop} \label {lifting}
If $(\gamma ,\kappa )$ is a dual pair of orbits in $Z$ with
$z_0\in \gamma \cap \kappa $ and $\tilde z_0$ is in a closed
$G_{z_0}$-orbit in the fiber $F:=\pi ^{-1}(z_0)$, then
$\tilde \gamma :=G.z_0$ and $\tilde \kappa :=K^\mathbb C.z_0$
define a dual pair $(\tilde \gamma ,\tilde \kappa )$ in 
$\tilde Z$. Furthermore, the mapping 
$\pi \vert \tilde \gamma :\tilde \gamma \to \gamma $ is
proper.
\end {prop}
\begin {proof}
Since $K_{z_0}$ is a maximal compact subgroup of $G_{z_0}$, it
acts transitively on the compact orbit $G_{z_0}.\tilde z_0$.
Consequently $\tilde \gamma \cap F=K_{z_0}.\tilde z_0$ and, since
$\gamma \cap \kappa =K.z_0$, it follows that
$\tilde \gamma \cap \tilde \kappa =K.\tilde z_0$.

\s
The properness of $\pi \vert \tilde \gamma $ follows immediately
from the fact that it is a homogenous fibration with compact
fiber.
\end {proof}

In the following result we maintain the above notation.

\begin {prop} 
If $\gamma $ is not closed, then
$C(\gamma )\subset C(\tilde \gamma )$
\end {prop}

\begin {proof}
Let $\{ g_n\} $ be a sequence in $C(\tilde \gamma )$ which converges
to $g\in \text{bd}(C(\tilde \gamma )$.  By Prop.\ref {real cycles} there
is a sequence $\tilde z_n$ in $g_n(\tilde \kappa )\cap \tilde \gamma $
which diverges in $\tilde \gamma $.  Since
$\pi \vert \tilde \gamma $ is proper, the corresponding sequence
$z_n$ in $C(\gamma )$ is also divergent.  

\m
Now either $g_n\not \in C(\gamma )$ infinitely often, in which case
$g\not \in C(\gamma )$ or we may assume that $\{ g_n\}\subset C(\gamma )$.
In the latter case, since $z_n\in g_n(\kappa )\cap \gamma $ is divergent
in $\gamma $ it likewise follows that $g\not \in C(\gamma )$.

\s
Hence $\text{bd}(C(\tilde \gamma ))\cap C(\gamma )=\emptyset $ and
the desired result follows.
\end {proof}
 
Finally, let $\gamma $ be such that 
$\text{c}\ell (\gamma )$=$\gamma \ \dot \cup \ \gamma _{\text{c}\ell }$
and $\tilde \gamma $ be as above.
If $g_n \in C(\gamma )$ converges to $g\in \text{bd}(C(\gamma ))$, then,
again by Prop.\ref {real cycles}, it follows that 
$g(\text{c}\ell (\kappa ))\cap \gamma _{\text{c}\ell }\not =\emptyset $.  
Consequently, $C(\gamma _{\text{c}\ell })\subset C(\gamma )$.

\s
Thus, in the notation above, we have the following consquence.

\begin {cor} \label {closed orbits}
$C(\gamma _{\text{c}\ell })\subset C(\gamma )\subset C(\tilde \gamma )$  
\end {cor}

\subsection {Proof of the main theorem}\label {proof}

\m
Let us first consider the case where $\gamma $ is not closed and
let $S\in {\cal S}_\kappa $ as above. Since 
$\kappa \cap S \cap \gamma $ already realizes the intersection
number $S.\text{c}\ell(\kappa )$, it follows that
$g(\text{c}\ell (\kappa ))\cap Y=\emptyset $ for all
$g\in C(\gamma )$.

\s
From now on we regard the homogeneous space $\Omega =G^\mathbb C/K^\mathbb C$
as a family of $q$-dimensional cycles, $q:=dim_\mathbb C\kappa $,
via the identification of $g\in C(\gamma )$ with the cycle
$g(\text{c}\ell (\kappa ))$.

\m
Since $Y:=S\setminus O$ can be given the structure of a 
very ample divisor in $S$, it follows that incidence variety
of all $q$-dimensional cycles $C$ with $C\cap Y\not =\emptyset $ contains 
a complex hypersurface which pulls back to a complex hypersurface
$H\subset \Omega $ (\cite {BK},\cite {BM}). 

\m
By the method of (\cite {BK}), for $f\in \Gamma (S,{\cal O}(*Y))$ an
appropriately chosen meromorphic function with poles along $Y$,
the trace-transform $Tr(f)$ has a non-empty polar set $P$ which is
the union of a certain number of components of $H$. Apriori it is
possible that there are cycles $g\in C(\gamma )$ with
$\text{c}\ell (g(\kappa ))\cap Y\not =\emptyset $, but with 
$\vert \text{c}\ell (g(\kappa ))
\cap \Sigma \vert =S.\text{c}\ell (g(\kappa ))$.  These
would necessarily have positive-dimensional intersection with
$S$ and the non-discrete components would lie in $S\setminus \Sigma $.

\s
However, this phenomenon occurs on a lower-dimensional subvariety
in $\Omega $, and such cycles are limits of generic cycles which
intersect $S$ only in points of $\Sigma $ and of course at only
finitely many points which are bounded away from $Y$.  
Thus the value of $Tr(f)$ at 
$g(\text{c}\ell (g(\kappa ))$ is finite and therefore such
cycles are not in $P$. As a consequence 
$P\cap C(\gamma )=\emptyset $ and we redefine $H$ to be $P$ in
the sequel.

\begin {prop} \label {inclusions}
If $\gamma $ is not closed, 
then there exists a Borel subgroup $B$ which contains a factor $AN$ of
an Iwasawa-decomposition $G=KAN$ and a $B$-invariant hypersurface
$H$ in $\Omega $ so that 
$\Omega _{AG}\subset C(\gamma )\subset \Omega _H $
\end {prop}
\begin {proof}
From the characterization of cycle spaces of the open $G$-orbits
(\cite {HW1},\cite {FH}) it follows in particular that the intersection
of all such cycle spaces for the open orbits in $G^\mathbb C/B$
is $\Omega _{AG}$ (This intersection property also follows
in most cases from the results in \cite {GM}).  By Prop. 8.1 of (\cite {GM}) 
it then follows that

\begin {equation*}
\Omega _{AG}\subset C(\gamma ).
\end {equation*} 

In the case of a non-closed orbit we have
the $B$-invariant hypersurface $H$ in the complement of $C(\gamma )$
and, since $C(\gamma )$ is $G$-invariant and contains the
base point in $\Omega _{AG}$ it follows that it   
is contained in the connected component of
\begin {equation*}
\Omega \setminus (\underset{k\in K}{\cup}k(H))
\end {equation*}
which by definition is $\Omega _H$.
\end {proof}

\begin {thm} \label {equality}
If $G$ is not of Hermitian type, then
\begin {equation*}
C(\gamma )=\Omega _{AG}
\end {equation*}
for all $\gamma \in Orb_Z(G)$.
\end {thm}
\begin {proof}
Since $G$ is non-Hermitian, $\Omega _H$ is Kobayashi hyperbolic
(see \cite {H},\cite {FH}) and the main theorem of (\cite {FH}) can
therefore be applied: A $G$-invariant, Stein, Kobayashi hyperbolic
domain containing $\Omega _{AG}$ is $\Omega _{AG}$ itself.
Thus $\Omega _{AG}=\Omega _H$ and the result for non-closed orbits
follows from from Prop.\ref {inclusions}.

\m
By Cor. \ref {closed orbits} it then follows that
$C(\gamma _{\text{c}\ell })\subset \Omega _{AG}$.  But,
if $\tilde Z=G^\mathbb C/B$, it is clearly the case
that $C(\tilde \gamma _{\text{c}\ell })\subset C(\gamma _{\text{c}\ell })$
and it is known that $C(\tilde \gamma _{\text{c}\ell })=\Omega _{AG}$
(\cite {B},\cite {H}, see also \cite {FH}). Thus,
$C(\gamma _{\text{c}\ell })\supset \Omega _{AG}$ and the result is proved for
closed orbits as well.
\end {proof}

\m
It should be remarked that we actually only use the Kobayashi
hyperbolicity of $\Omega _H$ and that one might expect the
following result in the Hermitian case: If $\Omega _H$ is not
Kobayashi hyperbolic, then $C(\gamma )$ can be identified with
the bounded Hermitian symmetric space associated with $G$.

\m
{\bf Note added in proof.}  With a few added details which we note here, 
we also handle the Hermitian case.  For this assume that $G$ is of
Hermitian type and let $P$ be a complex subgroup of $G^\mathbb C$
which properly contains $K^\mathbb C$.  Then $X:=G^\mathbb C/P$ is
a compact Hermitian symmetric space.  There are only two choices
for $P$ and the (open) $G$-orbit of the neutral point in $X$ is the
bounded Hermitian symmetric domain ${\cal B}$ or its conjugate
$\bar {\cal B}$.

\begin {prop} 
If $G$ is of Hermitian type, then either 
\begin {equation*}
C(\gamma )=\Omega _{AG}={\cal B}\times \bar {\cal B}
\end {equation*}
or the base dual cycle $\text{c}\ell (\kappa )$ is $P$-invariant
and $C(\gamma )$ is either ${\cal B}$ or $\bar {\cal B}$, depending
on the choice of sign.
\end {prop}
\begin {proof}
In the Hermitian case it is known that $\Omega _{AG}$ agrees with
${\cal B}\times \bar {\cal B}$ in its natural embedding in
$G^\mathbb C/K^\mathbb C$ (\cite {BHH}).

\s
Let $\pi :G^\mathbb C/K^\mathbb C\to G^\mathbb C/P=:X$, where $P$ is
one of the two choices mentioned above, and let $B$ be an
Iwasawa-Borel subgroup of $G^\mathbb C$ which is being used for
the incidence geometry.  If the $B$-invariant hypersurface
$H$ in $G^\mathbb C/K^\mathbb C$ of $\S\ref {proof}$
is not a lift $H=\pi ^{-1}(H_0)$, then $\Omega _H$ is Kobayashi 
hyperbolic and therefore $\Omega _H=\Omega _{AG}$ (\cite {FH}). 
Thus, just as in the non-Hermitian case,
the desired result follows from Prop.$\,$\ref {inclusions}.

\m
In order to complete the proof, we assume that $H$ is a lift, but
that $\text{c}\ell (\kappa )$ is not $P$-invariant, and reach
a contradiction.

\m
Let $x_0\in \gamma $ be a base point with $\kappa =K^\mathbb C.x_0$.
Since $\kappa $ is not $P$-invariant, $\text{c}\ell (P.x_0)$
contains $\text{c}\ell (\kappa )$ as a proper algebraic subvariety.
Now the intersection $\kappa \cap \gamma $ is transversal in $Z$.
Thus every component of $P.x_0\cap {\cal O}$ is positive-dimensional.
Since ${\cal O}=\mathbb C^{m({\cal O})}$ is affine, every such
component has at least one point of $Y$ in its closure.

\m
Thus for $B$ an arbitrarily small neighborhood of the identity
in $G^\mathbb C$ there exists $p\in P$ and $g\in B$ with
$gp.x_0\in Y$, and consequently $gpK^\mathbb C\in \Omega _H$.

\m
But $H$ was assumed to be a lift.  This is equivalent to $\Omega _H$
being a lift, i.e., at the group level $\Omega _H$ is invariant 
under the action of $P$ defined by right-multiplication.  Hence
$gK^\mathbb C$ is also in $\Omega _H$.  However, $g$ can be chosen
arbitrarily near the identity, contrary to $C(\gamma )$ being
an open neighborhood of the neutral point in $G^\mathbb C/K^\mathbb C$.
\end {proof} 

\m
{\bf Acknowledgement.} We would like to thank the referee for his 
constructive remarks; in particular, for underlining the fact that the
set of $g\in G^\mathbb C$ with $g(\kappa )\cap \gamma $ non-empty
and compact could theoretically be rather bad, e.g., not open.
\hfill $\square $

{\begin{tabular}{ll}
AH and BN: \\
Fakult\" at f\" ur Mathematik\\
Ruhr--Universit\" at Bochum \\
D-44780 Bochum, Germany \\
                               &                                      \\
{\tt ahuck@cplx.ruhr-uni-bochum.de}\\
{\tt ntatin@cplx.ruhr-uni-bochum.de}
\end{tabular}}

\end {document}